\newcommand{\arxiv}[2]{#1} %for arxiv

\arxiv{\documentclass[a4paper,10pt]{article}}{\documentclass{mjm}}

\usepackage{paper-en}
\arxiv{\usepackage{hyperref}}{\usepackage{mjm}}

\def\thetitle{Tubed embeddings}
\def\theauthors{Anton Petrunin}

\arxiv{\hypersetup{colorlinks=true,
citecolor=black,
linkcolor=black,
anchorcolor=black,
filecolor=black,
menucolor=black,
urlcolor=black,
pdftitle={\thetitle},
pdfauthor={\theauthors}
}}{}

\begin{document}

\arxiv{\title{\thetitle}
\author{\theauthors}
\date{}
\maketitle}{\title{\thetitle}
\headtitle{\thetitle}
\headauthor{\theauthors}
\maketitle}

\begin{abstract}
We consider the following question:
When does a Riemannian manifold admit an embedding with a uniformly thick tubular neighborhood in another Riemannian manifold of large dimension?
\end{abstract}

\section*{Introduction}

\paragraph{}\label{par:Observation}
Let $M$ and $N$ be Riemannian manifolds; we will always assume that Riemannian manifolds are complete connected, and $C^\infty$-smooth.

A smooth embedding $f\:M\hookrightarrow N$ will be called \emph{tubed} if the image $f(M)\z\subset N$ is a closed subset and it admits a tubular $\eps$-neighborhood for some $\eps>0$;
that is, the closest-point projection from the $\eps$-neighborhood of $f(M)$ to $f(M)$ is uniquely defined.
Equivalently, one may say that \textit{$f(M)$ is a set with positive reach} \cite{federer,bangert}.

We say that $M$ has \emph{uniformly polynomial growth} if there is a polynomial $p$ such that 
\[\vol B(x,R)_M\le p(R)\]
for any point $x\in M$ and any radius $R\ge 0$;
here $B(x,R)_M$ denotes open ball in $M$ of radius $R$ centered at $x$, and $\vol$ stands for volume.
In this case, we say that $p$ is a \emph{growth polynomial} of $M$.

\begin{thm}{Observation}
Suppose a complete Riemannian manifold $M$ admits an isometric tubed embedding in a Euclidean space.
Then $M$ has
\begin{enumerate}[(i)]
\item bounded sectional curvature,
\item positive injectivity radius, and
\item\label{uniformly polynomial growth} uniformly  polynomial growth.
\end{enumerate}
\end{thm}

For example,
the hyperbolic plane has exponential growth;
therefore, it does not admit an isometric tubed embedding into $\RR^d$ for any integer $d$.

\paragraph{}\label{par:main} 
The sectional curvature of $M$ will be denoted by $\sec_M$;
its injectivity radius at point $p\in M$ will be denoted by $\inj_pM$, and $\inj M$ will denote the injectivity radius of $M$; that is, 
\[\inj M\df \inf\set{\inj_pM}{p\in M}.\]

We say that $M$ has \emph{bounded geometry} if $\inj M>0$ and for every integer $k\ge 0$ there is a constant $C_k > 0$ such that $|\nabla^k \Rm| \le C_k$;
here $\nabla$ denotes the Levi-Civita connection, and $\Rm$ is the curvature tensor on $M$.
(Evidently, manifolds with bounded geometry have bounded sectional curvature.
Some authors define bounded geometry as bounded sectional curvature and injectivity radius --- we do not follow this convention.)

\begin{thm}{Main theorem}
Suppose $M$ is a Riemannian manifold with bounded geometry.
Then, $M$ has uniformly  polynomial growth
if and only if there is an isometric tubed embedding $M \hookrightarrow \RR^d$ for some integer $d$.

Moreover, the dimension $d$ can be found in terms of the dimension and a growth polynomial of $M$.
\end{thm}

\paragraph{Other ambient spaces.}\label{par:other-intro}
Loosely speaking, the following statement says that there is no universal space --- no space can serve as a target of an isometric tubed embedding for \textit{any} $n$-manifold with bounded geometry.

\begin{thm}{Proposition}
Let $N$ be a Riemannian manifold with bounded geometry and $n\ge 2$.
Then there is an $n$-dimensional  Riemannian manifold $M$ with bounded geometry that does not admit an isometric tubed embedding in $N$.

Moreover, if $N_1, N_2, \dots $ is a sequence of Riemannian manifolds with bounded geometry of arbitrary dimensions.
Then there is an $n$-dimensional Riemannian manifold $M$ with bounded geometry that does not admit an isometric tubed embedding in any $N_k$.
\end{thm}

Note that in the second statement the dimensions of $N_i$ can be arbitrary.
In particular, there are manifolds with bounded geometry that do not admit tubed embedding into any hyperbolic space, or any product of hyperbolic planes and so on; see also §\ref{par:Hn+H2n}.

\paragraph{About the proofs.}
The observation in §\ref{par:Observation} is proved in §§\ref{lem:gauss}--\ref{par:obs-proof}.
It is done by means of elementary differential geometry.

To prove the main theorem (§\ref{par:main}), we approximate the manifold by a graph (§\ref{par:intersection-graph}) and apply a slight modification of 
the theorem by Robert Krauthgamer and James Lee \cite{krauthgamer-lee0,krauthgamer-lee1} (§\ref{par:Phi}).
It produces a map $f_1\:M\to \RR^{d_1}$ such that $|f_1(x)-f_1(y)|$ is bounded away from zero for all pairs $x,y\in M$ on distance at least $R$ from each other for some fixed $R>0$ (§\ref{par:rough-embedding}).

Further, in §\ref{par:local-embedding}, we construct another map $f_2\:M\to \RR^{d_2}$ that is bi-Lipschitz in every ball of radius $R$.
The direct sum $f_1\oplus f_2$ is an embedding.
We can choose $\eps>0$ such that $\eps^2\cdot h<\tfrac12\cdot g$,
where $h$ is the induced Riemannian metric on $M$ by $f_1\oplus f_2$ and $g$ is the original Riemannian metric on $M$.
Applying Nash's embedding theorem for metric $g-\eps^2\cdot h$,
we find another map $f_3\:M\to \RR^{d_3}$ such that the direct sum $f=(\eps\cdot f_1)\oplus (\eps\cdot  f_2)\oplus f_3$ is isometric.

We show that the constructions of $f_1$, $f_2$, and $f_3$ produce \textit{uniformly smooth} maps; see  §\ref{par:rough-embedding}.
In particular, the normal curvatures of $f$ are bounded.
Then we use elementary differential geometry to show that $f$ is tubed. 

The proof of the proposition in §\ref{par:other-intro} is a slight modification of the argument of Florian Lehner \cite[1.2]{lehner}; see also \cite{462670}. 
It states that \textit{no connected graph with bounded degree contains a copy of any connected graph of degree 3}; see §§\ref{par:univeral}--\ref{par:other-proof}.

\section*{Observation}

\paragraph{Lemma.}\label{lem:gauss}
\textit{Suppose $M$ is a smooth submanifold in $\RR^d$ with normal curvatures at most 1.
Then $-2\le \sec_M(\sigma)\le 1$ for any sectional direction $\sigma$ of $M$.}
\medskip

\parit{Proof.}
Denote by $\T_p$ and $\N_p$ the tangent and normal spaces at $p\in M$.
Denote by $s$ the second fundamental form of $M$ at $p$;
it is a symmetric quadratic form on $\T_p$ with values in $\N_p$.

By the assumptions, $|s(v,v)|\le 1$ for any $|v|=1$.
We need to show that curvature of $M$ at $p$ in any sectional direction
lies in the range $[-2,1]$.
Passing to a 2-dimensional subspace of $\T_p$, we may assume that $\dim M=2$.

Let $h\in \N_p$ and $\zh\in\RR$ be the average values of $s(u,u)$ and $|s(u,u)|^2$ for unit vectors $u\in \T_p$.
Since $\dim M=2$, the Gauss formula \cite[2.1]{petrunin2023} can be written as
\[K=3\cdot |h|^2-2\cdot \zh;\]
here $K$ is the Gauss curvature of $M$ at $p$.
Note that $|h|^2\le \zh$.
Since the normal curvatures are at most 1, we have $\zh\le 1$.
Hence the lemma follows.
\qeds

It is straightforward to check that the inequality $0\le |h|^2 \le \zh$ is the only restriction on $|h|^2$ and $\zh$.
Therefore, \textit{the obtained bounds are optimal}.

\paragraph{}\label{par:obs-proof}
The first two parts of the observation in §\ref{par:Observation} follow from a more general result by Stephanie Alexander and Richard Bishop \cite[1.4]{alexander-bishop},
but we include proofs for completeness.

\parit{Proof of the observation in §\ref{par:Observation}.}
Suppose  a smooth $n$-dimensional submanifold $M$ in $\RR^d$ has a $\eps$-thick tubular neighborhood $W$.

Note that normal curvatures of $M$ cannot exceed $\tfrac1\eps$.
By §\ref{lem:gauss}, 
\[-\tfrac2{\eps^2}\le\sec_M\le\tfrac1{\eps^2}.\]

\begin{wrapfigure}{o}{25 mm}
\vskip-6mm
\centering
\includegraphics{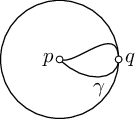}
\end{wrapfigure}

Let us show that $\inj M\ge\eps$; suppose $\inj_pM<\eps$ for some~$p$.
Since $\sec_M\le\tfrac1{\eps^2}$,
Klingenberg's lemma \cite[5.6]{cheeger-ebin} implies that there is a geodesic loop $\gamma$ based at $p$ with length $2\cdot \inj_pM$.
Let \[r=\max\set{|p-q|_{\RR^d}}{q\in\gamma}.\]
Observe that $r<\inj_pM$, and $\gamma$ has curvature at least $\tfrac1r$ at $q$.
It follows that the normal curvature of $M$ at $q$ in the direction of $\gamma$ exceeds $\tfrac1\eps$ --- a contradiction.

Consider the closest point projection $\pi\:W\z\to M$.
By Weyl's tube formula \cite{gray}, there is a constant $c\z=c(\eps,n,d)>0$ such that
\[\vol_d (\pi^{-1}X)\ge c\cdot\vol_n X\]
for any open set $X\subset M$.

{

\begin{wrapfigure}{o}{45 mm}
\vskip-14mm
\centering
\includegraphics{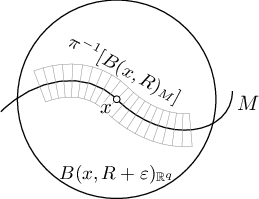}
\vskip5mm
\end{wrapfigure}

Note that $B(x,R)_M\subset B(x,R)_{\RR^d}$ for any $x\in M$;
here $B(x,R)_M$ and $B(x,R)_{\RR^d}$ denote balls in the metric of $M$ and $\RR^d$ respectively.
Therefore 
\[\pi^{-1}[B(x,R)_M]\subset B(x,R+\eps)_{\RR^d}.\arxiv{}{\qquad\qquad\qquad\qquad\qquad\qquad}\]
Since $R\mapsto \vol_d B(x,R+\eps)_{\RR^d}$ is a polynomial function, we conclude that $M$ has uniformly  polynomial growth.
\qeds

}

\section*{Graph embedding}

\paragraph{}\label{par:graph-embedding}
Let $\Gamma$ be a connected graph.
Recall that $\Gamma$ comes with the shortest-path distance on its vertex set $\Vert \Gamma$.
(In other words, we equip all graphs with the length metric having all edges of length 1, but consider only the restriction of this metric to the set of vertices and disregard the edges.)
So, we can talk about the ball $B(x,R)_\Gamma$ for any center $x\in \Vert \Gamma$ and radius $R>0$.
The number of vertices in $B(x,R)_\Gamma$ will be denoted by $|B(x,R)_\Gamma|$.
We say that $\Gamma$ has \emph{uniformly polynomial growth} if there is a polynomial $p$ such that 
\[|B(x,R)_\Gamma|\le p(R)\]
for any point $x\in \Gamma$ and any radius $R\ge 0$.
In this case, we say that $p$ is a \emph{growth polynomial} of $\Gamma$.
Note that \textit{any graph with uniformly polynomial growth has bounded degree}.

We will write $v\sim w$ if two vertices $v,w\in \Vert \Gamma$ are adjacent.
Let us denote by $\ZZ^d_\infty$ the graph with the vertex set formed by integer lattice $\ZZ^d\subset \RR^d$ and
\[v\sim w \iff  \|v-w\|_\infty=1.\]
The next statement follows from \cite[Theorem 5.5]{krauthgamer-lee1}.

\begin{thm}{Graph-embedding theorem}\label{thm:graph-embedding}
Suppose $\Gamma$ is a connected graph with uniformly polynomial growth.
Then, $\Gamma$ is isomorphic to a subgraph of $\ZZ^d_\infty$ for some integer~$d$.
Moreover, the dimension $d$ can be found in terms of a growth polynomial of~$\Gamma$.
\end{thm}

\paragraph{}\label{par:Phi}
Let us denote by $\Conv X$ the convex hull of the subset $X\subset \RR^d$.

\begin{thm}{Corollary}\label{cor:graph-embedding}
Suppose $\Gamma$ is a connected graph with uniformly polynomial growth.
Then for some $\rho>0$ and integer $d>0$ there is a map $\Phi\:\Vert \Gamma\to \RR^d$ such that
\begin{enumerate}[(a)]
\item if $v\sim w$, then $|\Phi(v)-\Phi(w)|\le \rho$, and 
\item if $V$ and $W$ are disjoint cliques in $\Gamma$, then the minimal distance from $\Conv[\Phi(V)]$ to $\Conv[\Phi(W)]$ is at least $1$.
\end{enumerate}
Moreover, the dimension $d$ and the distance $\rho$ can be found in terms of a growth polynomial of~$\Gamma$.
\end{thm}

\parit{Proof.}
By §\ref{par:graph-embedding}, it is sufficient to prove the corollary for $\Gamma=\ZZ^n_\infty$.

Take $d=n+2^n$.
Let us color vertices of $\ZZ^n_\infty$ in $2^n$ colors so that adjacent vertices get different colors.
Let $e_1,\dots,s_{2^n}$ be the standard basis in $\RR^{2^n}$.
Consider the map $\Phi\:\ZZ^n\to \RR^d=\RR^n\oplus \RR^{2^n}$ defined by $\Phi\:v\mapsto v\oplus e_{i(v)}$ where $i(v)$ denotes the color of $v$.
After an appropriate rescaling, the obtained map meets the conditions.
\qeds

\section*{Approximation by graph}

\paragraph{Intersection graph.}\label{par:intersection-graph}
Let $M$ be a complete connected Riemannian manifold.
Choose $r>0$.
Suppose $V\subset M$ is a maximal set of points (with respect to inclusion) at a distance at least $r$ from each other.
Given $v\in V$, set $B_v=B(v,r)_M$ and $\lambda\cdot B_v=B(v,\lambda\cdot r)_M$.
By the construction, we have that 
\textit{
\begin{enumerate}[(i)]
\item\label{item:1/2B} the balls $\{\tfrac12\cdot B_v\}_{v\in V}$ are disjoint;
\item $\{B_v\}_{v\in V}$ is an open cover of $M$;
\item\label{item:Lebesgue} $\{2\cdot B_v\}_{v\in V}$ is an open cover of $M$ with a Lebesgue number~$r$.
\item\label{item:1/N} If $M$ has bounded geometry, then for any positive real $\lambda$, there is an integer $N_\lambda$ such that any $\lambda\cdot r$-ball in $M$ contains at most $N_\lambda$ points of $V$.
\end{enumerate}
}

Let $\Gamma_\lambda$ be the intersection graph of the covering $\{\lambda\cdot B_v\}$ for $\lambda\ge 1$;
that is, $\Vert \Gamma_\lambda=V$, and 
\[v\sim w\quad\iff\quad \lambda\cdot B_v\cap \lambda\cdot B_w \ne \emptyset;\]
here $v\sim w$ means that $v$ is adjacent to $w$ in $\Gamma_\lambda$.

Assume $M$ has bounded geometry.
By \textit{\ref{item:1/N}}, $\Gamma_\lambda$ has degree at most~$N_{2\cdot\lambda}-1$.
Since $v\sim w$ implies |$v-w|_M< 2\cdot\lambda\cdot r$, we have
\[2\cdot\lambda\cdot r\cdot |v-w|_{\Gamma_\lambda}\ge |v-w|_M\eqlbl{eq:dist-G}\]
for any $v,w\in V$.

If $M$ has bounded geometry, there is $\delta>0$ such that $\vol (\tfrac12\cdot B_v) >\delta$ for any $v\in V$.
Therefore, \textit{\ref{item:1/2B}} and \ref{eq:dist-G} imply the following.

\textit{
\begin{enumerate}[(i)]\addtocounter{enumi}{4}
\item If $M$ has bounded geometry and uniformly polynomial growth, then the graph $\Gamma_\lambda$ has uniformly polynomial growth for any $\lambda\ge1$.\\
Moreover, a growth polynomial of $\Gamma_\lambda$ can be found in terms of $\lambda$, $r$ and a growth polynomial of $M$.
Namely, if $x\mapsto p(x)$ is a growth polynomial of $M$, then $x\mapsto p(2\cdot\lambda\cdot r\cdot x)$ is a growth polynomial of $\Gamma_\lambda$.
\end{enumerate}
} 

\section*{Factors of embedding}

\paragraph{}\label{par:rough-embedding}
Let $M$ be a smooth Riemannian manifold.
A map $f\:M\to\RR^d$ is called \emph{uniformly smooth} if 
for any integer $k\ge 1$ there is a constant $C_k$ such that  
\[|\tfrac{d^k}{dt^k}(f\circ\gamma)(t)|\le C_k\] 
for all unit-speed geodesics $\gamma$.

\begin{thm}{Proposition}
Let $M$ be a complete $n$-dimensional Riemannian manifold with bounded geometry and uniformly  polynomial growth.
Suppose $|\sec_M|\le \tfrac1{100}$ and $\inj M\ge 10$.
Then for some integer $d_1$ there is a uniformly smooth map $f_1\:M\z\to\RR^{d_1}$ such that 
\[|x-y|_M\ge 1
\quad\Longrightarrow\quad
|f_1(x)-f_1(y)|_{\RR^{d_1}}\ge 1.\]
Moreover, the dimension $d_1$ can be found in terms of a growth polynomial of~$M$.
\end{thm}

\parit{Proof.}
Let us apply the construction in §\ref{par:intersection-graph} for $r=\tfrac14$;
we obtain a subset $V\subset M$, balls $\{\lambda\cdot B_v\}_{v\in V}$, integer values $N_\lambda$ and connected graphs $\Gamma_\lambda$ with $\Vert\Gamma_\lambda =V$ for all $\lambda\ge 1$.

Now, let us apply the construction of partition of unity.
Choose a $C^\infty$-smooth function $\sigma\:\RR\to [0,1]$
such that $\sigma(t)=1$ for $t\le 1$ and $\sigma(t)=0$ for $t\ge 2$.
Further, let $\psi_v=\sigma(\tfrac{\dist_v}{r})$, for each $v\in V$;
note that the functions $\phi_v$ are uniformly smooth.
Set 
\[\Psi=\sum_{v\in V} \psi_v\quad \text{and}\quad \phi_v=\frac{\psi_v}{\Psi}.\]
The obtained functions $\phi_v$ form a partition of unity subordinate to $\{2\cdot B_v\}$.
By \textit{\ref{item:1/N}} in §\ref{par:intersection-graph}, $1\z\le \Psi\z\le N_2$; therefore, the functions $\phi_v$ are uniformly smooth. 

Given $x\in M$, consider the set 
\[V_x=\set{v\in V}{\phi_v(x)>0}.\]
If $\phi_v(x)>0$, then $2\cdot B_v\ni x$.
Therefore $V_x$ is a clique in $\Gamma_2$.
Moreover,
\[|x-y|\ge 4\cdot r=1
\qquad\Longrightarrow\qquad
V_x\cap V_y=\emptyset.
\eqlbl{eq:cliques}
\]

Let $\Phi\: V\to\RR^{n}$ be the map constructed in §\ref{par:Phi}.
Given $x\in M$, define $f_1(x)$ as the barycenter of points $\Psi(v)$ with masses $\phi_v(x)$ for each $v\in V_x$.
By §\ref{par:Phi}, $\diam V_x\le \rho$,
and since the functions $\phi_v$ are uniformly smooth, so is $f_1$.

By the construction, $f_1(x)$ lies in the convex hull of $\Phi(V_x)$ for any $x\in M$.
It remains to apply \ref{eq:cliques} and §\ref{par:Phi}.
\qeds

\paragraph{Proposition.}\label{par:local-embedding}
\textit{
Let $M$ be a complete $n$-dimensional Riemannian manifold with bounded geometry.
Suppose $|\sec_M|\le \tfrac1{100}$ and $\inj M\ge 10$.
Then for some integer $d_2$, there is a uniformly smooth map $f_2\:M\to\RR^{d_2}$ such that 
$f_2$ is uniformly bi-Lipschitz in all unit balls.
}

\parit{Proof.}
Choose $r=1$ and apply the constructions in §\ref{par:intersection-graph}.

Let $\sigma\:\RR\to[0,1]$ be as in §\ref{par:rough-embedding}.
For each $v\in V$, choose an isometry $\iota_v\:\T_v\z\to\RR^n$.
Consider the map $s_v\:M\to \RR^n$ defined by
\[s_v\:\exp_vx\mapsto \sigma(2\cdot |x|)\cdot \iota_v(x).\]

Since $\Gamma_2$ has degree at most $N_2-1$,
the set $V$ can be subdivided into $N_2$ subsets $V_1,\dots, V_{N_2}$ such that 
for any $i$, the balls $2\cdot B_v$ for $v\in V_i$ are disjoint.
Consider the map $S_i\:M\to\RR^n$ defined by
\[S_i\df\sum_{v\in V_i}s_v.\]
Note that $S_i$ is uniformly smooth and it is uniformly bi-Lipschitz in every ball $2\cdot B_v$ for $v\in V_i$.
It follows that $S_1\oplus\dots\oplus S_{N_2}\:M\to \RR^{N_2\cdot n}$ is uniformly smooth and bi-Lipschitz in  $2\cdot B_v$ for any $v\in V$.
The proposition follows since $r=1$ is a Lebesgue number of the covering $\{2\cdot B_v\}$; see \textit{\ref{item:Lebesgue}} in §\ref{par:intersection-graph}.
\qeds

\section*{Uniform Nash's construction}

\paragraph{Uniformly regular metrics.}
\label{par:uniformly-regular}
Let $M$ be a smooth $n$-dimensional manifold with a chosen locally finite atlas $\mathcal{A}$.
We say that $\mathcal{A}$ is \emph{appropriate} if it meets the following conditions:
\begin{itemize}
\item Each chart $\sigma\in \mathcal{A}$ is given by a map from an open set $U_\sigma\subset M$ to the unit ball $\DD^n\subset \RR^n$.
\item The covering $\{U_\sigma\}$ of $M$ has bounded multiplicity.
\item There exists $r\in (0,1)$ such that $\{\sigma^{-1}(r\cdot \DD^n)\}_{\sigma\in\mathcal{A}}$ is a cover for $M$.
\item All the transition maps of $\mathcal{A}$ are uniformly bi-Lipschitz and for any integer $k\ge 0$ they  have uniformly bounded $C^k$-norms. 
\end{itemize}

A metric tensor $g$ on $M$ will be called \emph{uniformly regular} with respect to an appropriate atlas $\mathcal{A}$
if for every integer $k\ge 0$, all components of the metric tensor in all charts of $\mathcal{A}$ have bounded $C^k$-norms,
and all the charts are uniformly bi-Lipschitz with respect to the metric induced by $g$ on $M$. 

Analogously, we can define a \emph{uniformly regular sequence of metrics} on a fixed manifold with a chosen appropriate atlas.
(Note that a smooth metric on a fixed compact manifold, say $\mathbb{S}^n$, is uniformly smooth for some appropriate atlas,
but a sequence of smooth metrics does not have to be uniformly smooth.)

Let $M$ be a Riemannian manifold with bounded geometry.
Applying rescaling if necessary, we may assume $|\sec_M|\le \tfrac1{100}$ and $\inj M\ge 10$.

Let us apply the construction in §\ref{par:intersection-graph} for $r=\tfrac12$.
Given $v\in V$, consider the normal chart $\sigma_v\:2\cdot B_v\to\DD^n$ with the center at $v$;
we will call such an atlas \emph{standard}.
It is straightforward to check that this atlas is appropriate and
the metric tensor is uniformly regular with respect to this atlas.
In particular, \textit{any manifold of bounded geometry has a uniformly regular metric in an appropriate atlas}.
The following statement from \cite{disconzi-shao-simonett} implies that the converse holds as well.

\begin{thm}{Proposition}
A Riemannian manifold $M$ has bounded geometry if and only if its metric is uniformly regular for some appropriate atlas.
\end{thm}

Also note that \textit{if $M$ comes with a uniformly regular metric for some appropriate atlas $\mathcal{A}$,
then a map $f\:M\to \RR^d$ is uniformly smooth if and only if for every integer $k\ge 0$,
the functions $f\circ\sigma-f\circ\sigma(0)\: \DD^n\to \RR^d$ for charts $\sigma\in \mathcal{A}$
have uniformly bounded $C^k$-norms.}

\paragraph{Proposition.}\label{par:nash}
\textit{Let $M$ be a complete $n$-dimensional Riemannian manifold with bounded geometry.
Then for some integer $d_3$, there is a uniformly smooth isometric immersion $f_3\:M\to\RR^{d_3}$.}

\medskip

The statement follows from the proof of Nash's theorem about regular embedding;
we will describe the needed changes in \cite{nash}.

Suppose we have a sequence of Riemannian metrics $g_1,g_2,\dots $ on the sphere $\mathbb{S}^n$.
By Nash's theorem, there are smooth embeddings $s_1,s_2,\z\dots\:\mathbb{S}^n\to \RR^d$
such that each $s_i$ induces $g_i$ on $\mathbb{S}^n$; here $d\z=\tfrac12\cdot n\cdot(3\cdot n+11)$.

Suppose that the sequence of metrics $g_i$ is uniformly regular with respect to a finite appropriate atlas on $\mathbb{S}^n$.
Then the construction in \cite[Part C]{nash} with the same choices of functions produces uniformly smooth embeddings $s_i$.
Namely, the functions $\psi^r$ in \cite[(C2)]{nash} and $\alpha_r$ in \cite[(C10)]{nash} are uniformly smooth and we can assume that $\lambda$ in \cite[(C12)]{nash} is the same for all $g_i$. 

Further, the so-called \emph{Nash's reduction} \cite[Part D]{nash} reduces the existence of an isometric immersion of $M$ to the existence of isometric immersions $(\mathbb{S}^n,g_i)\to \RR^d$ for a sequence of Riemannian metrics $g_i$.
Nash starts with an atlas with bounded multiplicity, and for each chart, he produces a sequence of metrics $g_i$ on $\mathbb{S}^n$ with a sequence of smooth maps $\phi_i\:M\to (\mathbb{S}^n,g_i)$ 
such that the sum of metrics induced by $\phi_i$ on $M$ coincides with $g$.
Then an isometric immersion is constructed from embeddings $s_i$ as above.

Applying rescaling, we may assume that $|\sec_M|\le \tfrac1{100}$ and $\inj M\ge 10$.
Let us choose a standard atlas on $M$ as in §\ref{par:uniformly-regular}.
Then Nash's construction produces a uniformly smooth sequence of maps $\phi_i$ and a sequence of metrics satisfying the above condition.
As a result, we get a uniformly smooth isometric immersion $(M,g)\to \RR^{d_3}$, where $d_3= \tfrac12\cdot N_2\cdot n\cdot(3\cdot n+11)$.

\section*{Proof assembling}

\paragraph{Proof of the main theorem (§\ref{par:main}).}\label{par:main-proof}
Let $f_1\:M\to\RR^{d_1}$ and $f_2\:M\to\RR^{d_2}$ be the maps constructed in §\ref{par:rough-embedding} and §\ref{par:local-embedding}.
Consider the map $f_1\oplus f_2\:M\z\to \RR^{d_1+d_2}$ defined by $p\mapsto (f_1(p),f_2(p))\in \RR^{d_1}\oplus\RR^{d_2}=\RR^{d_1+d_2}$.

Let $h$ be the metric induced by $f_1\oplus f_2$
and $g$ be the original metric on $M$.

Choose $\eps>0$ such that $\eps^2\cdot h<\tfrac12\cdot g$.
The metric tensor $g-\eps^2\cdot h$ is uniformly regular and smooth on $M$ for a standard atlas; see §\ref{par:uniformly-regular}.
By §\ref{par:uniformly-regular}, $(M,g-\eps^2\cdot h)$ has bounded geometry.
By §\ref{par:nash}, there is a uniformly smooth map $f_3\:M\to\RR^{d_3}$ that induces metric $g-\eps^2\cdot h$ on $M$.

Therefore, the map $f=(\eps\cdot f_1)\oplus (\eps\cdot f_2)\oplus f_3\:M\to \RR^{d_1+d_2+d_3}$ is a uniformly smooth isometric immersion.
In particular, the image $f(M)$ has bounded normal curvatures.
By §\ref{par:rough-embedding} and §\ref{par:local-embedding}, $f_1\oplus f_2$ is an embedding;
therefore so is $f$.
It remains to show that $f$ is tubed.

\begin{wrapfigure}{o}{45 mm}
\vskip-0mm
\centering
\includegraphics{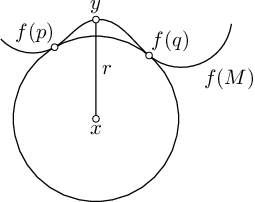}
\end{wrapfigure}

Assume the contrary;
that is, for any $\eps>0$ there is a point $x\in \RR^d$ that lies at a  distance smaller than $\eps$ to $f(M)$ such that the distance function from $x$ has at least two minimum points on $f(M)$,
say $f(p)$ and $f(q)$.
Denote by $\gamma$ a minimizing geodesic from $p$ to $q$ in $M$.

Suppose $\eps<\tfrac12$.
The properties in §\ref{par:rough-embedding} and §\ref{par:local-embedding} imply that 
\[\length\gamma<L\cdot \eps\]
for some fixed $L\in \RR$.
Let $y\in f(\gamma)$ be a point that maximizes the distance from $x$, and $r=|x-y|$.
The triangle inequality implies that $r\le \eps\cdot(1+L)$.
It follows that the curvature of $f(\gamma)$ at $y$ is at least $\tfrac1r$.
The latter is impossible if $\eps$ is small.
\qeds

\section*{Universal space}

\paragraph{Theorem.}\label{par:univeral}\textit{
Let $\Gamma$ be a connected graph of bounded degree.
Then there is a connected graph $\Delta$ of degree $\le 3$ such that there is no map $f\:\Vert \Delta\to \Vert \Gamma$ such that for some constant $K$ we have the following.
\begin{enumerate}[(i)]
\item $|f^{-1}(h)|\le K$ for any vertex $h$ in $\Gamma$.
\item $|f(v)-f(x)|_\Gamma\le K$ for any two adjacent vertices $v, w$ in $\Delta$.
\end{enumerate}
}

\textit{Moreover, if $\Gamma_1,\Gamma_2,\dots$ be a sequence of connected graphs with bounded degrees.
Then there is a connected graph $\Delta$ of degree $\le 3$ such that for any $i$ there is no map $f\:\Vert \Delta\to \Vert \Gamma_i$ with above properties.}

\medskip

A slightly weaker statement was proved by Florian Lehner \cite[1.2]{lehner} and rediscovered by Fedor Petrov \cite{462670};
we will modify their proof straightforwardly.
The conditions in the theorem are closely related to the so-called \emph{regular maps};
see \cite{hume-mackay-tessera} and the references therein.

\parit{Proof.}
Suppose $\Gamma$ is an infinite connected graph with bounded degree.
Observe that $\Gamma$ has a countable set of vertices; so we can identify $\Vert \Gamma$ with the set of natural numbers $\NN$.

Choose a sequence $n_1,n_2,\dots$ of natural numbers such that each $n\in \NN$ appears infinitely many times in it;
for example, we may take the following sequence: $1, 2, 1, 3, 1, 2, 1, 4, 1,\dots$

Let us construct a graph $\Delta$ with degree at most $3$ such that no map $f$ meets the condition in the theorem for any $K$.
In the following construction, the set of vertices of $\Delta$ will be $\NN$, and it will contain an infinite path $P=1\sim 2\sim \dots$
There will be a fast growing sequence $S_1<S_2<\dots$ such that 
if $i\sim j$ for some $i<j$, then $j-i=1$ or $S_k<i<j\le S_{k+1}$ for some $k$.

The subgraphs induced by $\{1,\dots, S_k\}$ will be denoted by $\Delta_k$;
let $P_k$ be the corresponding path $1\sim 2\sim \dots \sim S_k$.
The graphs $\Delta_k$ will be constructed recursively so that they meet the following condition:
\begin{clm}{}\label{clm:fk}
There is no map $f_k\:\Vert \Delta_k\to \Vert \Gamma$ that meets the conditions for $K=k$ and such that $f_k(1)\z=n_k$.
\end{clm}
Once it is done, the theorem follows. 

We assume that $\Delta_0$ is the empty graph, and $S_0=0$.
Assume $\Delta_{k-1}$ is constructed, let us  construct $\Delta_k$.
To do so, we need to choose large $S_k\z\gg S_{k-1}$ and add diagonals to the path $S_{k-1}+1\sim\dots\sim S_k$.

Suppose $L<(S_k-S_{k-1})/2$.
Let $A$ ($B$) be the first (respectively, the last) $L$ elements in the integer interval from $S_{k-1}+1$ to $S_k$.
We may choose an arbitrary bijection $A\leftrightarrow B$ and connect corresponding vertices by edges.
This way we get $L!$ different graphs; in particular, the total number of graphs $\Delta_k$ has superexponential growth with respect to $S_k$.
More precisely, we count labeled graphs; that is, graphs $\Delta_k$ with vertices labeled by $\{1,\dots,S_k\}$ and a path $P_k=1\sim 2\sim\dots\sim S_k$.

How much of these graphs have the required maps $f\:\Vert \Delta_k\to \Vert \Gamma$?
Suppose that the degree of $\Gamma$ does not exceed $d-1$.
Then each $k$-ball in $\Gamma$ has at most $d^k$ vertices.
Therefore we have at most $d^{k\cdot S_{k+1}}$ different maps of the path $P_k$.
Given one such map, there are at most $k\cdot d^k$ ways to add a diagonal at each vertex.
Therefore the number of labeled graphs that admit a required map $f\:\Vert \Delta_k\to \Vert \Gamma$ for $K=k$ grows exponentially with respect to $S_k$.
It follows that we may choose a sufficiently large $S_k$,
so that there is a graph $\Delta_k$ that meets \ref{clm:fk}.

To prove the second statement, note that the set of vertices in all $\Gamma_i$ is countable,
and therefore we can apply the same argument.
\qeds

\paragraph{Proof of the proposition in §\ref{par:other-intro}.}\label{par:other-proof}
Let $N$ be a manifold with bounded geometry and $n\ge 2$;
we need to construct an $n$-dimensional manifold $M$ with bounded geometry that does not admit an isometric tubed embedding into $N$.

Let $\Gamma_2$ be the intersection graph obtained by applying the construction in §\ref{par:intersection-graph} to $N$.
Recall that $\Gamma_2$ has bounded degree.
Let $\Delta$ be the graph provided by §\ref{par:univeral} for $\Gamma_2$;
it is a connected graph with  degree at most~3.

Let us choose a \emph{tube}; that is, the cylinder $\mathbb{S}^{n-1}\times [0,1]$ with a smooth Riemannian metric such that a neighborhood of each boundary component is isometric to a small annulus in the unit sphere $\mathbb{S}^n$, so it can be used to connect spheres.
Prepare a unit sphere $\mathbb{S}^n$ for every vertex in $\Delta$.
Connect two spheres by a tube if the corresponding vertices are adjacent.
The obtained manifold $M$ has bounded geometry.

Assume there is an isometric tubed embedding $\iota\colon M\hookrightarrow N$.
We can assume that $\Delta$ is embedded in $M$ so that each vertex lies on the corresponding sphere.
Consider the map $\Delta\to \Gamma_2$ that sends a vertex $v\in \Vert \Delta$ to a vertex of $\Gamma_2$ that lies at minimal distance from $\iota(v)$.
The obtained map satisfies the conditions in §\ref{par:univeral} for some $K$ --- a contradiction.

To prove the second statement, argue in the same manner, applying the second statement in §\ref{par:univeral}.
\qeds

\section*{Final remarks}

\paragraph{}\label{par:remarks}
Let us say that a Riemannian manifold $M$ has \emph{$k$-bounded geometry} if it has positive injectivity radius and there is a constant $C_k > 0$ such that $|\nabla^i \Rm| \le C_k$ for any nonnegative integer $i\le k$.
For example, $0$-bounded geometry means positive injectivity radius and bounded curvature.

In the main theorem, it is sufficient to assume that $M$ has only $3$-bounded geometry;
the proof requires only minor changes.
Below we discuss $1$- and $0$-bounded geometry.

Let us generalize the definitions of appropriate atlas and uniformly regular metric (see §\ref{par:uniformly-regular}).
Choose an integer $m\ge 1$ and $\alpha\in(0,1)$.
Suppose that in the definition of appropriate atlas,
instead of bounded $C^k$ norm of the transition maps for all $k$ 
we only require bounded $C^{m+1,\alpha}$ norm.
Let us call such atlas \emph{$C^{m+1,\alpha}$-appropriate}.
Furthermore, in the definition of uniformly regular metric let us assume that the atlas is only $C^{m+1,\alpha}$-appropriate, and instead of bounded $C^k$ norm for all $k$,
we only require that metric has bounded $C^{m,\alpha}$ norm in every chart.
Let us call such metrics \emph{uniformly $C^{m,\alpha}$-regular}.
Note that a $C^\infty$-smooth metric on a noncompact manifold might fail to be uniformly $C^{m,\alpha}$-regular even if injectivity radius is positive.

Recall that if metric is $C^{m,\alpha}$, then it remains to be $C^{m,\alpha}$ in the harmonic charts.
This fundamental result was obtained by Idzhad Sabitov and Samuil Shefel \cite{sabitov-shefel} and rediscovered by Dennis DeTurck and Jerry Kazdan \cite{deturck-kazdan}.
Applying this technique, one could get that \textit{manifolds with $1$-bounded geometry have uniformly $C^{2,\alpha}$-regular metric for any $\alpha\in (0,1)$}.

Arguing as in in §\ref{par:nash} with Günther's version of the embedding theorem \cite{guenther},
we obtain the following: \textit{for any $m\ge 2$ and $\alpha\in(0,1)$, any manifold with uniformly $C^{m,\alpha}$-regular metric admits a uniformly $C^{m,\alpha}$-smooth isometric immersion into Euclidean space;}
in particular, such an immersion has bounded normal curvatures.
Furthermore, following §\ref{par:other-proof}, we get that
\textit{for any $m\ge 2$ and $\alpha\in(0,1)$, any manifold with uniformly $C^{m,\alpha}$-regular metric and uniformly polynomial growth admits a uniformly $C^{m,\alpha}$-smooth tubed isometric embedding into a Euclidean space of large dimension.}

Straightforward calculations show that a uniformly $C^{1,1}$-regular metric has bounded sectional curvature and positive injectivity radius.
It is unknown if the converse holds.
Also it is unknown if a manifold with $C^{1,1}$-regular metric admits a $C^{1,1}$-smooth isometric embedding into Euclidean space.
These are two parts of the question posed by Vitali Kapovitch and Alexander Lytchak \cite[1.10]{kapovitch-lytchak}.
If both questions have affirmative answers, then the presented argument should produce a converse to the observation in §\ref{par:Observation}.
(It is not hard to check that if converse to the observation hods, then any compact Riemannian manifold with bounded sectional curvature admits a $C^{1,1}$-smooth isometric embedding into Euclidean space.)

If the smoothness is below $C^{1,1}$, then there is no hope for tubed isometric embedding.
However, it is easy to see that an embedding $f\:M\to \RR^d$ with bounded normal curvature is tubed if and only if it is \emph{uniformly graphical};
that is, there is $r>0$ such that for any $p\in M$ the set $B(f(p),r)\z\cap f(M)$ is a graph of a function $\RR^n\to \RR^{d-n}$ for a suitable choice of coordinate system in $\RR^d$.
Applying the result of Anders K\"{a}ll\'{e}n \cite{kallen}, the argument in §\ref{par:other-proof} produces that 
\textit{for any $\alpha\in(0,1)$, any manifold with $C^{1,\alpha}$-regular metric and uniformly polynomial growth admits a  uniformly $C^{1,\alpha}$-smooth uniformly graphical embedding into a Euclidean space of large dimension.}
In other words, for arbitrary small positive $\eps$, the $C^{2,\eps}$ and $C^{1,1-\eps}$ versions of the uniformly graphical question have affirmative answers, but $C^{1,1}$---the most important case---remains open.

\paragraph{Question}\label{par:Hn+H2n}
\textit{Suppose that an $n$-dimensional Riemannian manifold $M$ has bounded geometry.
What is necessary and sufficient condition for the existence of tubed isometric embeddings of $M$
into (1) hyperbolic space of sufficiently large dimension or (2) product of sufficiently many hyperbolic planes?}

\medskip

By §\ref{par:other-intro}, the conditions have to be nontrivial, 
but it says nothing about what these conditions should be.
The condition of uniformly polynomial growth is sufficient for both questions.
It follows from the main theorem; the second case is trivial and the first case follows 
since a horosphere defines a tubed isometric embedding of a Euclidean space into a hyperbolic space.
Evidently this condition is not necessary.

An argument presented by David Hume and Alessandro Sisto \cite[1.1]{hume-sisto} can be used to find some necessary conditions. For example, \textit{the product of the hyperbolic plane and the real line does not admit a tubed embedding into hyperbolic space of arbitrary dimension}.
Also, a manifold with infinite asymptotic dimension does not admit a tubed embedding into a product of hyperbolic planes.
The latter assertion relies on the monotonicity of asymptotic dimension with respect to maps between graphs satisfying the two conditions in §\ref{par:univeral}.
Namely, \textit{if $f\:\Delta\to \Gamma$ is such a map between graphs with bounded degrees, then $\mathrm{asdim}\,\Delta\le \mathrm{asdim}\,\Gamma$}.
(The proof is left to the reader.)

Finally, let us note that the proof in §\ref{par:other-proof} works equally well if $N$ has only 0-bounded geometry. 
It would be interesting to understand what happens if in §\ref{par:other-intro} we assume that $N$ is just a complete connected smooth Riemannian manifold. 

\paragraph{Acknowledgments.} 
This note is inspired by a question posed by Peter Michor \cite{124840}.
I want to thank 
\texttt{aorq} (\texttt{user:1079} at MathOverflow),
Alexander Dranishnikov,
Michael Gromov,
David Hume,
Nina Lebedeva,
James Lee,
Alexander Lytchak, 
Fedor Petrov,
and
Daniele Semola
for their help.

This work was partially supported by the National Science Foundation, grant DMS-2005279,
and the L. Euler International Mathematical Institute at PDMI RAS, grant 075-15-2022-289.

{\sloppy
\def\emph{\textit}
\printbibliography[heading=bibintoc]
\fussy
}
\end{document}